\documentclass{article}
\usepackage{arxiv}
\usepackage[colorlinks,citecolor=blue]{hyperref}   
\usepackage{url}            
\usepackage{microtype}      
\usepackage{cleveref}       
\usepackage{lipsum}         
\usepackage{caption}
\usepackage{cite}
\usepackage{doi}
\usepackage[T1]{fontenc}
\usepackage{authblk}

\usepackage{extarrows}		
\usepackage{tikz}			
\usepackage{esint}			

\begin{document}
	\title{On the $L_\mathrm{Y J}(\xi,\eta,X)$ constant \\for the Bana\'s-Fr\k{a}czek space}
	\author[a]{Yuxin Wang}
	\author[a]{Qi Liu\thanks{Qi Liu:liuq67@aqnu.edu.cn}}
	\author[a]{Linhui Chen}
	\author[a]{Xiewei Tan}
	\author[b]{Muhammad Sarfraz}
	\affil[a]{School of Mathematics and physics, Anqing Normal University, Anqing 246133,P.R.China}
	\affil[b]{School of Mathematics and System Sciences, Xinjiang University, Urumqi 830046, P.R.China}
	\date{}
	\maketitle
	
	\begin{abstract}
		In this paper, for any $\lambda\geq1,~R_\lambda^2$ is the Bana\'s-Fr\k{a}czek space. The
		exact value of $L_\mathrm{YJ}(\xi,\eta,X)$ for this space will be calculated. Specifically, $L_{\mathrm{YJ}}( \xi , \eta , R_{\lambda }^{2}) = 1+ \frac {2\xi \eta }{\xi ^{2}+ \eta ^{2}}( 1- \frac 1{\lambda ^{2}})$ is the result thereafter through meticilous computation.
	\end{abstract}
	
	\keywords{Banach spaces\and Geometric constants\and Bana\'s-Fr\k{a}czek space.}
	\textbf{Mathematics Subject Classification:}{ 46B20}
	
	\section{Introduction} 
	Throughout this article, let $X$ be a real Banach space with dimension at least 2. We will define $S_X$ as the unit sphere, that is $S_X=\{x\in X:\|x\|_X=1\};B_X$ as the closed unit ball, that is $B_X=\{x\in X:\|x\|_X\leq1\}.$ Additionally, we will use $ext(B_x)$ to indicate the collection of extreme points of $B_x.$\\The geometric properties of Banach spaces are studied in terms of their geometric constants. We can better comprehend and analyze the features of Banach spaces by using alternative geometric constants, which may give us distinct geometric characteristic information about the space. Investigating the exact value of geometric constants for particular spaces is therefore beneficial to us.
	
	The authors presented the space $R_{\lambda}^2$ in \cite{2}, which is considered to be $R_{\lambda}^{2}:=(R^{2},\|\cdot\|_{\lambda})$, where $\lambda>1$ and $$\|(a,b)\|=\max\left\{\lambda|a|,\sqrt{a^{2}+b^{2}}\right\}.$$ Furthermore, the authors named the space as the Bana\'s-Fr\k{a}czek space.

		The Bana\'s-Fr\k{a}czek space was generalized by C. Yang and X. Yang\cite{02} in 2016. $X_{\lambda,p}$ indicates this promotional space, which is thought to be where $\lambda>1$ and $p\geq1$, $X_{\lambda,p}:=(R^{2},\|\cdot\|_{\lambda,p})$, of which $\|(a,b)\|_{\lambda,p}=\max\left\{\lambda|a|,(a^{p}+b^{p})^\frac{1}{p}\right\}$, respectively.
	
	Let's now go over the definitions and findings of a few geometric constants
	that are relevant to the topic that follows.\\Remember that the definition of James constant $J(X) $\cite{12}  as shown below: $$J(X)=\sup\left\{\min\{\|x+y\|,\|x-y\|\}:x,y\in X,\|x\|=1,\|y\|=1\right\}.$$

	 With the deepening of the research, in \cite{1}, scholars have generalized the James constant $J(X)$, and the generalized form is as follows:$$J_{\lambda,\mu}(X)=\sup\Big\{\min\{\|\lambda x+\mu y\|,\|\mu x-\lambda y\|\}:x,y\in B_X\Big\}.$$
	
	Furthermore, in\cite{03}, the author defined the James type constant as$$J_{X,t}(\tau)=\sup\{\mu_t(\|x+\tau y\|,\|x-\tau y\|):x,y\in S_X\}.$$
	
	By careful calculation, C. Yang and X. Yang\cite{02} obtained the exact value of the James type constant for the $X_{\lambda,p}$ space. The results is as follows: 
	$$J_{X_{\lambda,p},t}(1)=\begin{cases}2^{1-\frac{1}{t}}\Big(1+\Big(1-\frac{1}{\lambda^{p}}\Big)^{\frac{t}{p}}\Big)^{\frac{1}{t}},&if~t\geq p,\\2^{1-\frac{1}{t}}\lambda\left(1+\lambda^{\frac{tp}{t-p}}\right)^{\frac{1}{t}-\frac{1}{p}},&ift<p~and~\lambda^{p}\leq1+\lambda^{\frac{p}{p}},\\2^{1-\frac{1}{t}}\Big(1+\Big(1-\frac{1}{\lambda^{p}}\Big)^{\frac{t}{p}}\Big)^{\frac{1}{t}},&if~t<p~and~\lambda^{p}\geq1+\lambda^{\frac{p}{p}}.\end{cases}$$
	
In \cite{01.1}, Clarkson  defined the von Neumann–Jordan constant of a Banach space X as the smallest constant C for which
   $$\frac{1}{C}\leq\frac{\|x+y\|^2+\|x-y\|^2}{2\|x\|^2+2\|y\|^2}\leq C,$$ where $x,y\in X$ and not both $0$.

Subsequently, the von Neumann-Jordan constant scholars are provided with an equivalent definition by study, as follows:

$$
C_{\mathrm{NJ}}(X) = \sup \left\{ \frac{\|x + y\|^2 + \|x - y\|^2}{2\|x\|^2 + 2\|y\|^2}:x,y\in X \text{~not both 0}\right\}.
$$
 
The von Neumann-Jordan constant $C_{\mathrm{NJ}}(X)$ has been shown to be a useful tool for characterizing Hilbert spaces, uniformly non-square spaces, and super-reflexive spaces after a great deal of investigation by scholars. For additional information and findings about $C_{\mathrm{NJ}}(X)$, see \cite{09,12}. Naturally, a lot of academics are also fixated on figuring out how much this constant is in particular spaces. The authors then gave the precise value of $C_\mathrm{NJ}(R^2_{\lambda})$ in \cite{04}, where $C_\mathrm{NJ}(R^2_{\lambda})=2-\frac{1}{\lambda^2}$.

While many scholars are studying the properties of the $C_{\mathrm{NJ}}(X)$  constant itself, other scholars have begun to generalize the constant. Thus, in \cite{06}, the generalized von Neumann-Jordan constant is defined as $$C_\mathrm{NJ}^{(p)}(X):=\sup\left\{\frac{\|x+y\|^p+\|x-y\|^p}{2^{p-1}\left(\|x\|^p+\|y\|^p\right)}:x,y\in X,(x,y)\neq(0,0)\right\},$$ where $1\leq p<\infty.$

The exact value of the $C_\mathrm{NJ}^{(p)}(X)$ constant for the Bana\'s-Fr\k{a}czek space, which was determined in 2018 by C. Yang and H. Li \cite{07}  was presented by where $\lambda>1,~p\geq2$, for any $\lambda$ such that $\lambda^2(1-\frac{1}{\lambda^2})^\frac{p}{2}\geq1$, $C_\mathrm{NJ}^{p}(R^2_\lambda)=1+(1-\frac{1}{\lambda^2})^\frac{p}{2}$.

	In a new research, Liu et al.\cite{08} introduced a new geometric constant with a skew connection that denotes to $L_\mathrm{Y J}(\xi,\eta,X)$ for $\xi,\eta>0$ as $$L_\mathrm{Y J}(\xi,\eta,X)=\sup \left\{\frac{\|\xi x+\eta y\|^2+\|\eta x-\xi y\|^2}{(\xi^2+\eta^2)(\|x\|^2+\|y\|^2)}: x,y\in X, (x,y)\neq (0,0) \right\}.$$
	
Additionally, the author \cite{11} provides a comparable constant $$L'_\mathrm{Y J}(\xi,\eta,X)=\sup \left\{\frac{\|\xi x+\eta y\|^2+\|\eta x-\xi y\|^2}{2(\|x\|^2+\|y\|^2)}: x,y\in S_X \right\}.$$

In recent years, a few scholars began to conduct research on the $L_\mathrm{Y J}(\xi, \eta, X)$ constant. The authors \cite{09} compared the relationship between the $L_\mathrm{Y J}(\xi, \eta, X)$ constant and the $J_{\lambda, \mu}(X)$ constant, and established an inequality between these two constants. Furthermore, in the same year, the authors \cite{10} obtain the exact value of the $L_\mathrm{Y J}(\xi, \eta, X)$ constant and the $L'_\mathrm{Y J}(\xi, \eta, X)$ constant for  the regular octagon space.

	Naturally we want to ask"what's $L_\mathrm{Y J}(\xi,\eta,X)$ constant for the Bana\'s-Fr\k{a}czek  space?" In this article, we calculate and prove that for any $\lambda>1$, $L_\mathrm{Y J}(\xi,\eta,R_{\lambda}^2)=1+\frac{2\xi \eta}{\xi^2+\eta^2}(1-\frac{1}{\lambda^{2}}).$ 
	
	\section{Main results}
	Before starting, we need to introduce the following lemma.
	\begin{Lemma}\cite{04}\label{l1}
		If $\lambda\geq\sqrt2$ and $|x_1|\leq\frac{1}{\lambda}$,$|y_1|\leq\frac{1}{\lambda}$,then$$(\lambda^2-1)|x_1y_1|+\sqrt{(1-x_1^2)}\sqrt{(1-y_1^2)}\leq2-\frac{2}{\lambda^2}.$$
	\end{Lemma}
	We will also need the following standard lemma.
	
	\begin{Lemma}\label{l2}
		Let $0\leq t\leq1$, $1\leq\lambda<\sqrt{2}$ and $0\leq x,y\leq\frac{1}{\lambda}$.\\(a) If $$f(x,y)=2t(\lambda^{2}-1)\xi\eta x y+2t\xi\eta\sqrt{1-x^{2}}\sqrt{1-y^{2}}+\lambda^{2}t^{2}\xi^{2}y^{2}+\lambda^{2}\eta^{2}x^{2},$$then 
		\begin{equation*}
			\max\{f(x,y):0\leq x,y\leq\frac{1}{\lambda}\}=f(\frac{1}{\lambda},\frac{1}{\lambda})=4t\xi\eta\cdot\frac{\lambda^2-1}{\lambda^2}+t^2\xi^2+\eta^2.
		\end{equation*} \\(b)If $$g(x,y)=2t(\lambda^{2}-1)\xi\eta x y+2t\xi\eta\sqrt{1-x^{2}}\sqrt{1-y^{2}}+\lambda^{2}t^{2}\eta^{2}y^{2}+\lambda^{2}\xi^{2}x^{2},$$then 
		\begin{equation*}
			\max\{g(x,y):0\leq x,y\leq\frac{1}{\lambda}\}=g(\frac{1}{\lambda},\frac{1}{\lambda})=4t\xi\eta\cdot\frac{\lambda^2-1}{\lambda^2}+t^2\eta^2+\xi^2.
		\end{equation*} 
	\end{Lemma}

	\begin{proof}
		Assume that $t\in[0,1]$. Suppose that $\max\{f(x,y):0\leq x,y\leq\frac{1}{\lambda}\}$ attains at some $(x,y)\in(0,\frac{1}{\lambda})\times(0,\frac{1}{\lambda})$. We have $$f_x=2\lambda^{2}\eta^2x+2t(\lambda^{2}-1)\xi\eta y-2tx\xi\eta\frac{\sqrt{1-y^{2}}}{\sqrt{1-x^{2}}}=0,$$ $$f_y=2\lambda^{2}\xi^{2}t^{2}y+2t(\lambda^{2}-1)\xi\eta x-2ty\xi\eta\frac{\sqrt{1-x^{2}}}{\sqrt{1-y^{2}}}=0,$$that is $$\begin{cases}\lambda^{2}\eta^{2}+t(\lambda^{2}-1)\xi n\frac{y}{x}=t\xi\eta\sqrt{\frac{1-y^{2}}{1-x^{2}}},\hspace{1.5em}(i)\\\lambda^{2}\xi^{2}t+(\lambda^{2}-1)\xi\eta\frac{x}{y}=\xi\eta\sqrt{\frac{1-x^{2}}{1-y^{2}}}.\hspace{1.8em}(ii)\end{cases}$$\\Now, through $(i)$ multiplied by $(ii)$, we have $$t\xi^{2}\eta^{2}=t\lambda^{4}\xi^{2}\eta^{2}+\lambda^{2}(\lambda^{2}-1)\xi^{3}n\frac{x}{y}+t^{2}\lambda^{2}(\lambda^{2}-1)\xi\eta^{3}\frac{y}{x}+t(\lambda^{2}-1)^{2}\xi^{2}\eta^{2},$$ and that is equivalent to $$0=2t\xi\eta \lambda^{2}+\xi^{2}\lambda^{2}\frac{x}{y}+t^{2}\lambda^{2}\eta^{2}\frac{y}{x},$$  which is a contradition.
		
		This means that $f(x,y)$ cannot attain its maximum value at $(x,y)\in(0,\frac{1}{\lambda})\times(0,\frac{1}{\lambda})$. Then we prove that  $\max\{f(x,y):0\leq x,y\leq\frac{1}{\lambda}\}=f(\frac{1}{\lambda},\frac{1}{\lambda})$. 
		
		According to the above, we know that the maximum value of $f(x,y)$  is either $f(0,0)$ ,$f(0,\frac{1}{\lambda})$ ,$f(\frac{1}{\lambda},0)$ or $f(\frac{1}{\lambda},\frac{1}{\lambda})$. Now we compare the magnitudes of these values.
		
		Since $$f(0,0)=2t\xi\eta,~~~  f(\frac1\lambda,\frac1\lambda)=4t\xi\eta\cdot\frac{\lambda^{2}-1}{\lambda}+t^{2}\xi^{2}+\eta^{2},$$
		$$f(0,\frac{1}{\lambda})=2t\xi\eta\sqrt{1-\frac{1}{\lambda^{2}}}+t^{2}\xi^{2},~~~f(\frac{1}{\lambda},0)=2t\xi\eta\sqrt{1-\frac{1}{\lambda^{2}}}+\eta^{2}.$$
		It's easy to know that $f(\frac{1}{\lambda},\frac{1}{\lambda})$ is the maximum value among these four values.
		Therefore,$$\max\{f(x,y):0\leq x,y\leq\frac{1}{\lambda}\}=f(\frac{1}{\lambda},\frac{1}{\lambda})=4t\xi\eta\cdot\frac{\lambda^2-1}{\lambda^2}+t^2\xi^2+\eta^2.$$\\The proof of (a) is completed.
		
		Assume that $t\in[0,1]$. Suppose that $\max\{g(x,y):0\leq x,y\leq\frac{1}{\lambda}\}$ attains at some $(x,y)\in(0,\frac{1}{\lambda})\times(0,\frac{1}{\lambda})$. We have $$g_x=2\lambda^{2}\xi^2x+2t(\lambda^{2}-1)\xi\eta y-2tx\xi\eta\frac{\sqrt{1-y^{2}}}{\sqrt{1-x^{2}}}=0,$$ $$g_y=2\lambda^{2}\eta^{2}t^{2}y+2t(\lambda^{2}-1)\xi\eta x-2ty\xi\eta\frac{\sqrt{1-x^{2}}}{\sqrt{1-y^{2}}}=0,$$that is $$\begin{cases}\lambda^{2}\xi^{2}+t(\lambda^{2}-1)\xi n\frac{y}{x}=t\xi\eta\sqrt{\frac{1-y^{2}}{1-x^{2}}},\hspace{1.5em}(iii)\\\lambda^{2}\eta^{2}t+(\lambda^{2}-1)\xi\eta\frac{x}{y}=\xi\eta\sqrt{\frac{1-x^{2}}{1-y^{2}}}.\hspace{1.8em}(iiii)\end{cases}$$\\Now, through $(iii)$ multiplied by $(iiii)$, we have $$t\xi^{2}\eta^{2}=t\lambda^{4}\xi^{2}\eta^{2}+\lambda^{2}(\lambda^{2}-1)\xi^{3}n\frac{x}{y}+t^{2}\lambda^{2}(\lambda^{2}-1)\xi\eta^{3}\frac{y}{x}+t(\lambda^{2}-1)^{2}\xi^{2}\eta^{2},$$ and that is equivalent to $$0=2t\xi\eta \lambda^{2}+\xi^{2}\lambda^{2}\frac{x}{y}+t^{2}\lambda^{2}\eta^{2}\frac{y}{x},$$  which is a contradition. 
		
		This means that $g(x,y)$ cannot attain its maximum value at $(x,y)\in(0,\frac{1}{\lambda})\times(0,\frac{1}{\lambda})$. Then we prove that  $\max\{g(x,y):0\leq x,y\leq\frac{1}{\lambda}\}=g(\frac{1}{\lambda},\frac{1}{\lambda})$. 
		
		According to the above, we know that the maximum value of $g(x,y)$  is either $g(0,0)$ ,$g(0,\frac{1}{\lambda})$ ,$g(\frac{1}{\lambda},0)$ or $g(\frac{1}{\lambda},\frac{1}{\lambda})$. Now we compare the magnitudes of these values.
		
		Since $$g(0,0)=2t\xi\eta,~~~  g(\frac1\lambda,\frac1\lambda)=4t\xi\eta\cdot\frac{\lambda^{2}-1}{\lambda}+t^{2}\eta^{2}+\xi^{2},$$
		$$g(0,\frac{1}{\lambda})=2t\xi\eta\sqrt{1-\frac{1}{\lambda^{2}}}+t^{2}\eta^{2},~~~g(\frac{1}{\lambda},0)=2t\xi\eta\sqrt{1-\frac{1}{\lambda^{2}}}+\xi^{2}.$$
		It's easy to know that $g(\frac{1}{\lambda},\frac{1}{\lambda})$ is the maximum value among these four values.
		
		Therefore,$$\max\{g(x,y):0\leq x,y\leq\frac{1}{\lambda}\}=g(\frac{1}{\lambda},\frac{1}{\lambda})=4t\xi\eta\cdot\frac{\lambda^2-1}{\lambda^2}+t^2\xi^2+\eta^2.$$\\The proof of (b) is completed.
	\end{proof}
	 
	 The following Theorem \ref{key} is our main result in this article.
	\begin{Theorem}\label{key}
		Let $\lambda\geq1$ {~and~} $R_{\lambda}^{2}$ is the Banas-Fraczek space. Then,\\$$L_\mathrm{Y J}(\xi,\eta,R_{\lambda}^{2})=1+\frac{2\xi \eta}{\xi^2+\eta^2}(1-\frac{1}{\lambda^{2}}).$$
	\end{Theorem}
	\begin{proof}To prove the theorem,  we prove that \begin{equation}\label{e1}
			{\frac{\|\xi x+\eta ty\|^2+\|\eta x-\xi ty\|^2}{(\xi^2+\eta^2)(1+t^2)}\leq 1+\frac{2\xi \eta}{\xi^2+\eta^2}(1-\frac{1}{\lambda^{2}})},
		\end{equation} holds for any $x,y\in ext(Bx)$ and any $t\in [0,1]$ first.
		
		Then we	consider where $\lambda\geq\sqrt{2}$ .
		
		Assume that $\lambda\geq\sqrt{2}$. Note that $ext(Bx)=\left\{(z_1,z_2):z_1^2+z_2^2=1,|z_1|\leq\frac{1}{\lambda}\right\}.$ \\To prove (\ref{e1}),we consider the following situations:\\\textbf{Case(1)} If
		$$
		\lambda\left|\xi x_{1}+\eta t y_{1}\right|\geq\sqrt{\left(\xi x_{1}+\eta t y_{1}\right)^{2}+\left(\xi x_{2}+\eta t y_{2}\right)^{2}},
		$$
		and
		$$
		\lambda\left|\eta x_{1}-\xi t y_{1}\right|\geq\sqrt{\left(\eta x_{1}-\xi t y_{1}\right)^{2}+\left(\eta x_{2}+\xi t y_{2}\right)^{2}}.
		$$
		Then we have
		\begin{equation*}
			\begin{split}
				\|\xi x+\eta t y\|^{2}+\|\eta x-\xi t y\|^{2}&=\lambda^{2}(\xi x_{1}+\eta t y_{1})^{2}+\lambda^{2}(\eta x_{1}-\xi ty_{1})^{2}\\&=\lambda^{2}(\xi ^{2}+\eta^{2})(x_{1}^{2}+t^{2}y_{1}^{2})\leq(\xi^{2}+\eta^{2})(1+t^{2}),
			\end{split} 
		\end{equation*}\\so \begin{equation}
			{\frac{\|\xi x+\eta ty\|^2+\|\eta x-\xi ty\|^2}{(\xi^2+\eta^2)(1+t^2)}\leq\frac{(\xi^2+\eta^2)(1+t^2)}{(\xi^2+\eta^2)(1+t^2)}\leq 1+\frac{2\xi \eta}{\xi^2+\eta^2}(1-\frac{1}{\lambda^{2}})}.
		\end{equation}
		\textbf{Case(2)} If
		$$
		\lambda\left|\xi x_{1}+\eta t y_{1}\right|\leq\sqrt{\left(\xi x_{1}+\eta t y_{1}\right)^{2}+\left(\xi x_{2}+\eta t y_{2}\right)^{2}},
		$$
		and
		$$
		\lambda\left|\eta x_{1}-\xi t y_{1}\right|\geq\sqrt{\left(\eta x_{1}-\xi t y_{1}\right)^{2}+\left(\eta x_{2}+\xi t y_{2}\right)^{2}}.
		$$   
		Then we have    
		\begin{equation*}
			\begin{aligned}||\xi x+\eta ty||^{2}+||\eta x-\xi ty||^{2}&=(\xi x_{1}+\eta ty_{1})^{2}+(\xi x_{2}+\eta ty_{2})^{2}+\lambda^{2}(\eta x_{1}-\xi ty_{1})^{2}\\&=\xi^{2}+\eta^{2}t^{2}+2t(\lambda^{2}-1)\xi \eta|x_{1}y_{1}|\\&+2t\xi \eta \sqrt{1-x_1^{2}}\sqrt{1-y_1^{2}}+\lambda^{2}t^{2}\xi^{2}y_1^{2}+\lambda^{2}\eta^{2}x_1^{2}.\end{aligned}
		\end{equation*}
		\textbf{Case(3)} If
		$$
		\lambda\left|\xi x_{1}+\eta t y_{1}\right|\geq\sqrt{\left(\xi x_{1}+\eta t y_{1}\right)^{2}+\left(\xi x_{2}+\eta t y_{2}\right)^{2}},
		$$
		and
		$$
		\lambda\left|\eta x_{1}-\xi t y_{1}\right|\leq\sqrt{\left(\eta x_{1}-\xi t y_{1}\right)^{2}+\left(\eta x_{2}+\xi t y_{2}\right)^{2}}.
		$$
		Then we have
		\begin{equation*}
			\begin{aligned}||\xi x+\eta ty||^{2}+||\eta x-\xi ty||^{2}&=\lambda^2(\xi x_{1}+\eta ty_{1})^{2}+(\eta x_{1}-\xi ty_1)^2(\eta x_{2}-\xi ty_{2})^{2}\\&=\xi^{2}t^2+\eta^{2}+2t(\lambda^{2}-1)\xi \eta|x_{1}y_{1}|\\&+2t\xi \eta \sqrt{1-x_1^{2}}\sqrt{1-y_1^{2}}+\lambda^{2}t^{2}\eta^{2}y_1^{2}+\lambda^{2}\xi^{2}x_1^{2}.\end{aligned}
		\end{equation*}
		For Case(2) and Case(3), by applying Lemma \ref{l1},we have \begin{equation}
			\begin{aligned}
				\frac{\|\xi x+\eta ty\|^2+\|\eta x-\xi ty\|^2}{(\xi^2+\eta^2)(1+t^2)}&\leq 1+\frac{2(2-\frac{2}{\lambda^2})\xi \eta t}{(\xi^2+\eta^2)(1+t^2)}\\&\leq 1+\frac{2\xi \eta}{\xi^2+\eta^2}(1-\frac{1}{\lambda^{2}}).
			\end{aligned}	
		\end{equation}
		\textbf{Case(4)} If
		$$
		\lambda\left|\xi x_{1}+\eta t y_{1}\right|\leq\sqrt{\left(\xi x_{1}+\eta t y_{1}\right)^{2}+\left(\xi x_{2}+\eta t y_{2}\right)^{2}},
		$$
		and
		$$
		\lambda\left|\eta x_{1}-\xi t y_{1}\right|\leq\sqrt{\left(\eta x_{1}-\xi t y_{1}\right)^{2}+\left(\eta x_{2}+\xi t y_{2}\right)^{2}}.
		$$
		Then we have
		\begin{equation*}
			\begin{aligned}\|\xi x+\eta y\|^{2}+\|\eta x-\xi y\|^{2}&=(\xi x_1+\eta ty_{1})^{2}+(\xi x_{2}+\eta ty_{2})^{2}\\&+(\eta x_{1}-\xi ty_{1})^{2}+(\eta x_{2}-\xi ty_{2})^{2}\\&=(\xi^{2}+\eta^{2})(1+t^{2}),\end{aligned}
		\end{equation*}\\so \begin{equation}\label{e4}
			{\frac{\|\xi x+\eta ty\|^2+\|\eta x-\xi ty\|^2}{(\xi^2+\eta^2)(1+t^2)}=\frac{(\xi^2+\eta^2)(1+t^2)}{(\xi^2+\eta^2)(1+t^2)}\leq 1+\frac{2\xi \eta}{\xi^2+\eta^2}(1-\frac{1}{\lambda^{2}})}.
		\end{equation}\\Then we prove that where $1\leq\lambda<\sqrt{2}$, (\ref{e1}) is also valid.
		
		Obviously, the above Case(1) and Case(4) are also valid for $1\leq\lambda<\sqrt{2}$. So we only need to consider Case(2) and Case(3) for $1\leq\lambda<\sqrt{2}$.\\In fact, by applying Lemma \ref{l2}, we have 	\begin{equation}\label{e5}
			\begin{aligned}||\xi x+\eta ty||^{2}+||\eta x-\xi ty||^{2}&=(\xi x_{1}+\eta ty_{1})^{2}+(\xi x_{2}+\eta ty_{2})^{2}+\lambda^{2}(\eta x_{1}-\xi ty_{1})^{2}\\&=\xi^{2}+\eta^{2}t^{2}+2t(\lambda^{2}-1)\xi \eta|x_{1}y_{1})|\\&+2t\xi \eta sin\sqrt{1-x^{2}}\sqrt{1-y^{2}}+\lambda^{2}t^{2}\xi^{2}y_1^{2}+\lambda^{2}\eta^{2}x_1^{2}\\&\leq\xi^{2}+\eta^{2}t^{2}+4t\xi\eta\cdot\frac{\lambda^2-1}{\lambda^2}+t^2\xi^2+\eta^2\\&=(\xi^{2}+\eta^{2})(1+t^{2}) +4t\xi\eta\cdot\frac{\lambda^2-1}{\lambda^2},\end{aligned}
		\end{equation}\\and\begin{equation}\label{e6}
			\begin{aligned}||\xi x+\eta ty||^{2}+||\eta x-\xi ty||^{2}&=\lambda^2(\xi x_{1}+\eta ty_{1})^{2}+(\eta x_{1}-\xi ty_1)^2(\eta x_{2}-\xi ty_{2})^{2}\\&=\xi^{2}t^2+\eta^{2}+2t(\lambda^{2}-1)\xi \eta|x_{1}y_{1})|\\&+2t\xi \eta sin\sqrt{1-x^{2}}\sqrt{1-y^{2}}+\lambda^{2}t^{2}\eta^{2}y_1^{2}+\lambda^{2}\xi^{2}x_1^{2}\\&\leq\xi^{2}t^2+\eta^{2}+4t\xi\eta\cdot\frac{\lambda^2-1}{\lambda^2}+\xi^2+\eta^2t^2\\&=(\xi^{2}+\eta^{2})(1+t^{2}) +4t\xi\eta\cdot\frac{\lambda^2-1}{\lambda^2}.\end{aligned}
		\end{equation}
		Combining(\ref{e5}) and(\ref{e6}), we have \begin{equation}\label{e7}
			\begin{aligned}
				\frac{\|\xi x+\eta ty\|^2+\|\eta x-\xi ty\|^2}{(\xi^2+\eta^2)(1+t^2)}&\leq 1+\frac{2(2-\frac{2}{\lambda^2})\xi \eta t}{(\xi^2+\eta^2)(1+t^2)}\\&\leq 1+\frac{2\xi \eta}{\xi^2+\eta^2}(1-\frac{1}{\lambda^{2}}).
			\end{aligned}	
		\end{equation}
		\hspace{1.5em}Therefore, combining (\ref{e4}) and (\ref{e7}), the conclusion of (\ref{e1}) is valid.
		
		On the other hand, if taking $x=\left(\frac{1}{\lambda},\sqrt{1-\frac{1}{\lambda^2}}\right)$,$y=\left(\frac{1}{\lambda},-\sqrt{1-\frac{1}{\lambda^2}}\right)$, we have $\|x\|=1$, $\|y\|=1$. Then\begin{equation}
			\frac{\|\xi x+\eta y\|^2+\|\eta x-\xi y\|^2}{2(\xi^2+\eta^2)}=1+\frac{2\xi \eta}{\xi^2+\eta^2}(1-\frac{1}{\lambda^{2}}).
		\end{equation}
		Obviously, if $\lambda\geq1$, we have $$|\xi +\eta|\geq\sqrt{\xi^2+\eta^2}\geq\sqrt{(\xi -\eta)^2+\frac{4\xi \eta}{\lambda^2}},$$ and $$\sqrt{(\xi +\eta)^2-\frac{4\xi \eta}{\lambda^2}}\geq\sqrt{(\xi+\eta)^2-4\xi \eta}=|\eta-\xi|.$$\\Hence, we have
		\begin{align*}
			\frac{\|\xi x+\eta y\|^2+\|\eta x-\xi y\|^2}{2(\xi^2+\eta^2)}=\frac{(\xi+\eta)^2+(\xi+\eta)^2-\frac{4\xi\eta}{\lambda^2}}{2(\xi^2+\eta^2)}&=\frac{2(\xi^2+\eta^2)+4\xi\eta-\frac{4\xi\eta}{\lambda^2}}{2(\xi^2+\eta^2)}&\\=1+\frac{2\xi \eta}{\xi^2+\eta^2}(1-\frac{1}{\lambda^{2}}).	
		\end{align*}
		\hspace{1.5em}Therefore, \begin{equation}\label{e9}
			L_\mathrm{Y J}(\xi,\eta,R_{\lambda}^2)\geq1+\frac{2\xi \eta}{\xi^2+\eta^2}(1-\frac{1}{\lambda^{2}}).
		\end{equation}
		By using (\ref{e1}) and (\ref{e9}), we complete the proof of the Theorem \ref{key}.	
	\end{proof}
			
In section 2, we have introduced that  $C_\mathrm{NJ}(R^2_{\lambda})=2-\frac{1}{\lambda^2}$ has given by Yang \cite{04}, in fact, the exact value of the $C_\mathrm{NJ}(X)$ constant for the Bana\'s Fr\k{a}czek space can also estimated by the above theorem, just as the following corollary states.
	\begin{Corollary}
		For any $\lambda\geq1$, $R^2_{\lambda}$ is the Bana\'s Fr\k{a}czek space. Then$$L_\mathrm{Y J}(1,1,R_{\lambda}^2)=C_\mathrm{NJ}(R^2_{\lambda})=2-\frac{1}{\lambda^2}.$$
	\end{Corollary}
	\begin{proof}
		Let $\xi=\eta=1$, then we have $L_\mathrm{Y J}(1,1,R_{\lambda}^2)=2-\frac{1}{\lambda^2}$ by Theorem \ref{key}. It's known that when $\xi=\eta=1$, the $L_\mathrm{Y J}(\xi, \eta, X)$ constant is just the $C_\mathrm{NJ}(X)$ constant, hence$$L_\mathrm{Y J}(1,1,R_{\lambda}^2)=C_\mathrm{NJ}(R^2_{\lambda})=2-\frac{1}{\lambda^2}.$$
	\end{proof}
	\begin{Corollary}
	Let $X$ be $R^2$  be furnished with the norm which is defined as $\|(a,b)\|=\max\left\{|a|,\sqrt{a^{2}+b^{2}}\right\}$. Then $$L_\mathrm{Y J}(\xi,\eta,X)=1.$$
	\end{Corollary}
	\begin{proof}
	By the preceding Theorem \ref{key}, let $\lambda=1$, we have $L_\mathrm{Y J}(\xi,\eta,R^2_{1})=1$. Then we obtain $L_\mathrm{Y J}(\xi,\eta,X)=1$ since  the Bana\'s-Fr\k{a}czek space refers to the aforementioned space $X$.
	\end{proof}
	
	In \cite{08}, the	author  provided a sufficient condition for judging the super-reflexivity of a Banach space through the $L_\mathrm{Y J}(\xi,\eta,X)$ constant. that is if $L_\mathrm{Y J}(\xi,\eta,X)<2,$ then  a Banach space $X$ is super-reflexive.\\From the relationship between the $L_\mathrm{Y J}(\xi,\eta,X)$ constant and super-reflexive spaces, we have the following corallary.
	\begin{Corollary}
		 The Bana\'s-Fr\k{a}czek space is super-reflexive where $\lambda\geq1$.
	\end{Corollary}
	\begin{proof}
According to Theorem \ref{key}, for $\lambda\geq1$,$$1+\frac{2\xi \eta}{\xi^2+\eta^2}(1-\frac{1}{\lambda^{2}})<2,$$is always valid. Hence the Bana\'s-Fr\k{a}czek space is super-reflexive.
	\end{proof}
	Furthermore, by considering whether the Bana\'s-Fr\k{a}czek space has weak normal structure, we draw the following corollary. Prior to presenting this corollary, we shall introduce a lemma provided by Liu.
	\begin{Lemma}\cite{08}\label{c}
		Let $X$ be a Banach space, if $L_\mathrm{Y J}(\xi, \eta, X)<\frac{(\xi+\eta)^2+(2\eta-\xi)^2}{2(\xi^2+\eta^2)}$ for $\eta\leq\xi\leq2\eta$, then $X$ has weak normal structure.
	\end{Lemma}
	\begin{Corollary}
		If $1\leq\lambda<\sqrt{\frac{4\xi\eta}{6\xi\eta-3\eta^2}}$, then for $\eta\leq\xi<\frac{3}{2}\eta$, the  Bana\'s-Fr\k{a}czek space has  weak normal structure.
		\end{Corollary}
	
		\begin{proof}
			For $\eta\leq\xi<\frac{3}{2}\eta$, if $1\leq\lambda<\sqrt{\frac{4\xi\eta}{6\xi\eta-3\eta^2}}$, then we have	$$L_\mathrm{Y J}(\xi,\eta,R_{\lambda}^2)=1+\frac{2\xi \eta}{\xi^2+\eta^2}(1-\frac{1}{\lambda^{2}})<\frac{(\xi+\eta)^2+(2\eta-\xi)^2}{2(\xi^2+\eta^2)}.$$
			This implies the  Bana\'s-Fr\k{a}czek space has  weak normal structure by Lemma \ref{c}.
		\end{proof}

\end{document}